\documentclass[11pt]{article}
\usepackage{amsmath, amssymb, latexsym, amsthm}

\newcommand{\sumlim}{\sum\limits}

\newtheorem{thm}{Theorem}[section]

\newtheorem{lem}[thm]{Lemma}
\newtheorem*{thmn}{Theorem}

\newtheorem{cor}[thm]{Corollary}

\newtheorem{obs}[thm]{Observation}

\newtheorem{exa}[thm]{Example}
\newtheorem{defn}[thm]{Definition}

\newtheorem{prop}[thm]{Proposition}
\newtheorem{conj}[thm]{Conjecture}
\newtheorem{clm}[thm]{Claim}
\newcommand{\een}{\end{enumerate}}
\newcommand{\blem}{\begin{lem}}
\newcommand{\elem}{\end{lem}}
\newcommand{\bcl}{\begin{clm}}
\newcommand{\ecl}{\end{clm}}
\newcommand{\bthm}{\begin{thm}}
\newcommand{\ethm}{\end{thm}}
\newcommand{\bpr}{\begin{prop}}
\newcommand{\epr}{\end{prop}}
\newcommand{\bco}{\begin{cor}}
\newcommand{\eco}{\end{cor}}
\newcommand{\bcon}{\begin{conj}}
\newcommand{\econ}{\end{conj}}
\newcommand{\bde}{\begin{defn}}
\newcommand{\ede}{\end{defn}}
\newcommand{\bex}{\begin{exa}}
\newcommand{\eexa}{\end{exa}}
\newcommand{\bobs}{\begin{obs}}
\newcommand{\eobs}{\end{obs}}
\newcommand{\bexe}{\begin{exe}}
\newcommand{\eexe}{\end{exe}}

\newcommand{\borel}{{\mathbb{B}^+}}
\newcommand{\borelc}{{\mathbb{B}_C^+}}

\renewcommand{\sp}{{Sp_{2n}(\mathbb{F})}}

\begin{document}

\title{Sign balance for finite groups of Lie type}

\author{Eli Bagno \\
\small Einstein institute of Mathematics\\
\small The Hebrew university\\
\small Givat Ram, Jerusalem, Israel. bagnoe@math.huji.ac.il\\
 Yona Cherniavsky \\
\small  Department of Mathematics and Statistics\\
 \small Bar-Ilan University\\
\small Ramat-Gan, Israel 52900 cherniy@macs.biu.ac.il}

\date{\today}

\maketitle

\begin{abstract}
 A product formula for the parity generating function of the number of $1'$s in invertible
 matrices over $\mathbb{Z}_2$ is given. The computation is based on algebraic tools such as the Bruhat decomposition.
 The same technique is used to obtain a parity generating function also
 for symplectic matrices over $\mathbb{Z}_2$. We present also a generating function for the sum of entries
 of matrices over an arbitrary finite field $\mathbb{F}_q$ calculated in $\mathbb{F}_q$. These formulas
are new appearances of the Mahonian distribution.

\end{abstract}



\bibliographystyle{is-alpha}

\section{Introduction}

Let $G$ be a subgroup of $GL_n(\mathbb{Z}_2)$. For every $K \in G$
define $o(K)$ to be the number of $1$'s in $K$. A natural problem
is to find the number of matrices with a given number of $1$'s, or
in other words, to compute the following generating function:
$$O(G,t)=\sumlim_{K \in G}{t^{o(K)}}.$$
In the case $G=GL_n(\mathbb{Z}_2)$, by considering the free action
of $S_n$ by permuting rows on the subset of $G$ containing the
matrices having a fixed number of $1'$s, it is not hard to see
that $O(G,t)$ has $n!$ as a factor but the complete generating
function can be rather hard to compute. A weaker variation of this
problem is to evaluate $O(G,-1)$. This is equivalent to
determining the difference between the numbers of even and odd
matrices, where a matrix is called {\it even} if it has an even
number of $1$'s and {\it odd} otherwise. The number $O(G,-1)$ will
be called the {\it parity difference} or the {\it imbalance} of
$G$. A group $G$ is called {\it sign-balanced} if $O(G,-1)=0$.

The notion of sign-balance has recently reappeared in a number of
contexts. Simion and Schmidt \cite{S} proved that the number of
$321$-avoiding even permutations is equal to the number of such
odd permutations if $n$ is even, and exceeds it by the Catalan
number $C_{\frac{1}{2}(n-1)}$ otherwise. Adin and Roichman
\cite{AR1} refined this result by taking into account the maximum
descent. In a recent paper \cite{St}, Stanley established the
importance of the sign-balance.\\

 In this paper we calculate the
parity difference for $G=GL_n(\mathbb{Z}_2)$. We also generalize
the problem of sign-balance to matrix groups over arbitrary finite
fields $\mathbb{F}_q$ where $q$ is a power of a prime $p$. It
turns out that the appropriate parameter for these fields is the
sum of nonzero entries of the matrix (mod $p$) rather than just
the
number of nonzero elements.\\

In the corresponding formulas there is an appearance of the number
$[n]_q=\frac{1-q^n}{1-q}=1+q+q^2+\cdots +q^{n-1}$. The most
important instance of this number is in the following theorem of
MacMahon:

\begin{prop}
$$\sum\limits_{\pi \in S_n}{q^{inv(\pi)}}=\sum\limits_{\pi \in S_n}{q^{maj(\pi)}}=[n]_q!$$
where $[n]_q!=[n]_q[n-1]_q \cdots [1]_q$, $inv(\pi)$ is the number
of inversions in $\pi$ (see definition in section \ref{type A}).
(The definition of $maj(\pi)$ is $maj(\pi)=\sumlim_{i=1}^n{i \cdot
\chi\{\pi(i+1)< \pi(i)\}}$.)

\end{prop}

A combinatorial parameter on permutations which has such a
distribution is called {\it Mahonian} after MacMahon.

Our results can also be seen as an example of the {\it cyclic
sieving phenomenon} (See \cite{RSW} for details). We  will expound
on this subject in the appendix.

We prove the following (See Theorems \ref{mainA} and \ref{mainC2})

\begin{thmn}
$$\sumlim_{K \in GL_n(\mathbb{Z}_2)}{(-1)^{o(K)}}=-2^{n \choose 2}[n-1]_2!$$
 \end{thmn}

\begin{thmn}
$$\sumlim_{K \in Sp_{2n}(\mathbb{Z}_2)}{(-1)^{o(K)}}=-2^{n^2}\cdot[2]_2[4]_2 \cdots [2n-2]_2$$
\end{thmn}

We also generalize the problem of sign-balance to matrix groups
over arbitrary finite fields. For $q$ a prime power the parameter
we work with is the sum of the entries of the matrix, calculated
in $\mathbb{F}_q$, rather than just the number of $1$-s. We order
the elements of $\mathbb{F}$, i.e., choose a bijection between
$\mathbb{F}$ and the set $\{0,...,q-1\}$ such that $0$ of
$\mathbb{F}$ goes to $0$. If we denote by $s(K)$ the sum of the
images of the entries of $K$ in $\mathbb{F}_q$ under this
bijection then the generating function we are interested in is
\begin{equation}\label{gen fun}
S(G,t)=\sum\limits_{K \in G}{t^{s(K)}} .
\end{equation}

For $G=GL_n(\mathbb{F}_q)$, instead of substituting $-1$ we
substitute the $q-th$ root of the unity in (\ref{gen fun}) to get:
(see theorem \ref{mainA}):

\begin{thmn}
$$\sumlim_{K \in GL_n(\mathbb{F}_q)}{\omega_q}^{s(K)}=(q-1)^{n-1}q^{n \choose 2}[n-1]_q!.$$
\end{thmn}

We note also that our results were obtained using algebraic tools
such as the Bruhat decomposition. This approach enables us to
present the generating function as a multiplicative formula which
turns out to have a Mahonian distribution, meanwhile giving a new
interpretation of the Mahonian distribution. We note that two other
approachs to the case of type A were proposed to us by Alex
Samorodnitzky and by the referee of the doctorate thesis of the first author.\\

The rest of this paper is organized as as follows: In Section
\ref{perl} we survey the Bruhat decomposition for type A and for
type C. In Section \ref{SignA} we prove our main theorem about the
sign balance for type A in the most generality i.e. for the groups
$GL_n(\mathbb{F}_q)$. In Section \ref{SignC} we present our
results for type C and in the appendix we extend about the
connection of our work to the Cyclic sieving phenomenon.

\section{Preliminaries} \label{perl}

\subsection{Finite groups of Lie type A}\label{type A}

 Let $\mathbb{F}$ be any field and let $G=GL_n(\mathbb{F})$ be
the group of invertible matrices over $F$. Let $H$ be the subgroup
of $G$ consisting of the diagonal matrices. This is a choice of a
torus for $G$. It is easy to show that the normalizer of $H$,
$N(H)$, is the group of  monomial matrices (i.e. each row and
column contains exactly one non-zero element). The quotient
$N(H)/H$ is called the Weyl group of type $A$, and is isomorphic
to $S_n$, the group of permutations on $n$ letters.
 The Borel subgroup $\mathbb{B}^+$ of the group $G$ consists of the
upper triangular matrices in $G$. The opposite Borel subgroup, consisting of
the lower triangular matrices, is denoted by $\mathbb{B}^-$. We denote by
$\mathbb{U}^+$ and $\mathbb{U}^-$ the groups of upper and lower triangular
matrices (respectively,) with $1$-s along the diagonal.

The Weyl group $S_n$ has a set
of Coxeter generators $\{s_1,...,s_{n-1}\}$, where $s_i$ can be
realized as the transposition $(i,i+1)$. We define also the length
of a permutation $\pi \in S_n$ with respect to the Coxeter generators to be:
$$\ell(\pi)=\min\{r \in N: \pi=s_{i_1} \cdots s_{i_r},
\mbox{for some } i_1,...,i_r \in [0,n-1]\}.$$

It is well known that for every $\pi \in S_n$
$$\ell(\pi)=inv(\pi)$$
where $$inv(\pi)=|\{(i,j)|\pi(i)>\pi(j),1 \leq i<j \leq n\}| .$$

\begin{prop}
For every finite field $\mathbb{F}$ with $q$ elements the order of $GL_{n}(\mathbb{F})$ is
$$q^{{n \choose 2}}(q-1)^{n}[n]_q!$$
\end{prop}

\subsection{The Bruhat Decomposition for type A}\label {BruhatA}

The Bruhat decomposition is a way to write every invertible matrix
as a product of two triangular matrices and a permutation matrix.
We start with the following definitions:

 For every permutation $\pi \in S_n$ we identify $\pi$ with the matrix:

\begin{center}
$$[\pi]_{i,j} =  \left\{ \begin{array}{cc}
{1}      &  { i=\pi(j)} \\
{0 }          & \text{otherwise}
\end{array}
\right. $$
\end{center}

Define for every $\pi \in S_n$:

$$\mathbb{U}_{\pi}=\mathbb{U}^- \cap (\pi \mathbb{U}^-\pi^{-1}).$$

$\mathbb{U}_{\pi}$ consists of the matrices with $1$-s along the
diagonal and zeros in place $(i,j)$ whenever $i<j$ or
$\pi^{-1}(i)<\pi^{-1}(j)$. This is an affine space of dimension
${n \choose 2}-\ell(\pi)$ over $\mathbb{F}$. ($\ell(\pi)$ is the
length of $\pi$ with respect to the Coxeter generators).


 Now, given $g \in G$, we can column reduce $g$ by
multiplying on the right by Borel matrices in order to get an
element $gb^{-1}$ satisfying the following condition:
\begin{equation}\label{representative} \tag{$\ast$}
\begin{split}
\text{The right most nonzero entry in each row is $1$} \\
\text{and it is the first nonzero entry in its column.}
\end{split}
\end{equation}

Those "leading entries" form a permutation matrix corresponding to
$\pi \in S_n$.

Now we can use $\pi^{-1}$ to rearrange the columns of $gb^{-1}$ in
order to get $gb^{-1}\pi^{-1}=u \in \mathbb{U}_{\pi}$, i.e., $g=u
\pi b$. This is called the {\it Bruhat decomposition} of the
matrix $g$. One can prove that this decomposition is unique, and
thus we have a partition of $G$ into double cosets indexed by the
elements of the Weyl group $S_n$.

If $\pi \in S_n$ then the double coset indexed by $\pi$ decomposes
into left $\mathbb{B}^+$-cosets in the following way: For every
choice of $u \in \mathbb{U}_{\pi}$, $u \pi$ is a representative of
the left coset $u\pi \borel$. Thus a general representative of the
double coset $\mathbb{U}_{\pi}$ can be taken as matrix of the form
(\ref{representative}), with every column filled with free
parameters beyond the leading $1$.

We summarize the information we gathered about the Bruhat
decomposition for type $A$ in the following:

\begin{prop}\label{summary-bruhat-A}
The group $GL_n(\mathbb{F})$ is a disjoint union of double cosets
of the form $\mathbb{U}_{\pi} \pi \borel$, where $\pi$ runs
through $S_n$. Every double coset decomposes into cosets of the
form $A\borel$ where $A$ is a general representative of the form
(\ref{representative}). The number of free parameters in $A$ is
equal to ${n \choose 2}-\ell(\pi)$.
\end{prop}

Here is an example of the coset decomposition for
$GL_3(\mathbb{Z}_2)$:

$$\mathbb{U}_1 1 \borel = \left\{\begin{pmatrix}
    1 & 0 & 0 \\
    \alpha & 1 & 0 \\
    \beta & \gamma & 1
    \end{pmatrix}\borel \mid \alpha,\beta,\gamma\in \mathbb{Z}_2\right\} $$ 
$$\mathbb{U}_{s_1} s_1 \borel = \left\{\begin{pmatrix}
    0 & 1 & 0 \\
    1 & 0 & 0\\
    \alpha & \beta & 1
    \end{pmatrix}\borel \mid \alpha,\beta \in \mathbb{Z}_2\right\} $$ 
$$ \mathbb{U}_{s_2} s_2\borel = \left\{\begin{pmatrix}
    1 & 0 & 0 \\
    \alpha & 0 & 1 \\
    \beta & 1 & 0
    \end{pmatrix}\borel \mid \alpha,\beta \in \mathbb{Z}_2\right\}  $$ 
$$\mathbb{U}_{s_2s_1} s_2s_1 \borel = \left\{\begin{pmatrix}
    0 & 1 & 0 \\
    0 & \alpha & 1 \\
    1 & 0 & 0
    \end{pmatrix}\borel \mid \alpha\in \mathbb{Z}_2\right\} $$ 
$$ \mathbb{U}_{s_1s_2} s_1s_2\borel = \left\{\begin{pmatrix}
    0 & 0 & 1 \\
    1 & 0 & 0 \\
    \alpha & 1 & 0
    \end{pmatrix}\borel \mid \alpha\in \mathbb{Z}_2\right\} $$ 
$$ \mathbb{U}_{s_1s_2s_1} s_1s_2s_1\borel  =  \begin{pmatrix}
    0 & 0 & 1 \\
    0 & 1 & 0 \\
    1 & 0 & 0
\end{pmatrix}\borel $$

\subsection{Lie Type C}\label{Background type C}

Let $J$ denote the $n \times n$ matrix
$$\begin{pmatrix}
0 & \cdot & \cdot & 1 \\
0 & \cdot & 1 & 0 \\
\cdot  & \cdot & \cdot & \cdot \\
1 & 0 & \cdot & 0
\end{pmatrix}$$

and let $$M=\begin{pmatrix}
0 & J \\
-J &  0
\end{pmatrix}.$$
The Lie group of type $C$, or the symplectic group, is defined
over the field $\mathbb{F}$ by:
$$Sp_{2n}(\mathbb{F})=\{A \in SL_{2n}(\mathbb{F}) \mid A^TMA=M\}.$$
This is the set of fixed points of the automorphism
$\varphi:SL_{2n}(\mathbb{F}) \longrightarrow SL_{2n}(\mathbb{F})$
given by: $\varphi(A)=M^{-1}(A^T)^{-1}M$.

An alternative way to present the symplectic group is the
following:  We define first a bilinear form on $\mathbb{F}^{2n}$:

\bde For every $x=(x_1,...,x_{2n}), y=(y_1,...,y_{2n}) \in
\mathbb{F}^{2n}$
$$B(x,y)=\sumlim_{i=1}^{n}{x_i \cdot y_{2n+1-i}}-\sumlim_{i=n+1}^{2n}{x_i \cdot y_{2n+1-i}}.$$
\ede

Denoting by $\{x_1,...,x_{2n}\}$ the set of columns of $X$ it is
easy to see that $X \in \sp$ if and only if the columns satisfy
the following set of equations:
\begin{center}
$$B(x_i,x_j) =  \left\{ \begin{array}{cc}
{(-1)^{i-j}}      &  {i+j=2n+1} \\
{0}          & {i+j\neq 2n+1}
\end{array}
\right.
$$
\end{center}

We end this section with the following well known fact:

\begin{prop}(See for example \cite[p.35]{O})

For every finite field $\mathbb{F}$ with $q$ elements the order of
$\sp$ is:
$$q^{n^2} (q-1)^{n}[2]_q \cdots [2n]_q$$
\end{prop}

\subsection{Bruhat Decomposition for Type C}\label{BruhatC}
In order to be able to present the Bruhat decomposition for type
C, we must first define a Borel subgroup for
$Sp_{2n}(\mathbb{F})$. We present this subject following
\cite{St}. Note that although the exposition of \cite{St} deals
with groups over algebraically closed fields, the results hold
also over finite fields. Start with the Borel subgroup
$\mathbb{B}^+$, chosen for type $A$, consisting of the upper
triangular matrices.

If $X=\left(\begin{smallmatrix}
A & B \\
0 &  C
\end{smallmatrix}\right) \in \mathbb{B}^+$, then $\varphi(X)=\left(\begin{smallmatrix}
J(C^T)^{-1}J & J(C^{T})^{-1}B^T(A^{T})^{-1}J \\
0 &  J(A^T)^{-1}J
\end{smallmatrix}\right) \in \mathbb{B}^+$.
(The automorphism $\varphi$ was defined in Section \ref{Background
type C}). Moreover, the automorphism $\varphi$ keeps the Borel
subgroup $\mathbb{B}^+$, as well as the groups of diagonal and
monomial matrices (denoted by $H$ and $N(H)$ respectively in
Section \ref{type A}) invariant. Thus we can take
$\mathbb{B}_C^+=Sp_{2n}(\mathbb{F}) \cap \mathbb{B}^+$ and
$\mathbb{B}_C^-=Sp_{2n}(\mathbb{F}) \cap \mathbb{B}^-$ as the
Borel subgroup and the opposite Borel subgroup of
$Sp_{2n}(\mathbb{F})$ respectively, and similarly for $H$ and
$N(H)$.

If we label the basis elements of the space on which $\sp$ acts by
indices \mbox{$n,n-1,...,1,-1,...,-n$,} then the Weyl group of
type $C$ can be realized as the group of those permutations $\pi
\in S_{2n}$ such that $\pi(-i)=-\pi(i)$. This is called the
octahedral group and denoted $B_n$. Looking at $B_n$ as a Coxeter
group, it has the following set of generators:
$$S=\{s_0,s_1,...,s_{n-1}\}$$ where $s_0$ is the transposition which permutes $1$ and $-1$ and
$s_i$ permutes $i$ and $i+1$, for $1 \leq i < n$. Just as in the
case of type $A$, we define here the length of a $B_n$-permutation
as:

$$\ell(\pi)=min\{r \in N: \pi=s_{i_1} \cdots s_{i_r},
\mbox{for some } i_1,...,i_r \in [0,n-1]\}.$$

We define also the groups $\mathbb{U}_C^+=\mathbb{U}^+ \cap
Sp_{2n}(\mathbb{F})$ and $\mathbb{U}_C^-=\mathbb{U}^- \cap
Sp_{2n}(\mathbb{F})$ to be the upper and lower unipotent subgroups
respectively. For every $\pi \in B_n$ we define
$\mathbb{U}^C_{\pi}=\mathbb{U}_C^- \cap (\pi \mathbb{U}_C^-
{\pi}^{-1})$. $\mathbb{U}^C_{\pi}$ is the intersection of $\sp$
with the set of matrices with $1$'s along the diagonal and zeros
at entries in location $(i,j)$ whenever $i<j$ or $\pi^{-1}(i) <
\pi^{-1}(j)$. This is an affine space of dimension $n^2
-\ell(\pi)$. (Here, $\ell(\pi)$ is the length function of $B_n$).

Now, we can use the Bruhat decomposition of $GL_{2n}(\mathbb{F})$
to produce the Bruhat decomposition for $Sp_{2n}(\mathbb{F})$. Let
$g \in Sp_{2n}(\mathbb{F})$. Consider $g$ as an element of
$GL_{2n}(\mathbb{F})$ and write $g=u \pi b$ where $\pi \in
S_{2n}$, $u \in \mathbb{U}_{\pi}$ and $b \in \borel$. We have:
$$g=\varphi(g)=\varphi(u)\varphi(\pi)\varphi(b),$$ but from the
uniqueness of the decomposition in $GL_{2n}(\mathbb{F})$ we have:
$$\varphi(u)=u,\quad  \varphi(\pi)= \pi h^{-1}, \quad
\varphi(b)=hb$$ where $h$ is diagonal and thus $\pi \in B_n$ and
$b \in \mathbb{B}_C^+$. This gives us the Bruhat decomposition.
The description of the double cosets and the coset representatives
is similar to the one given for type $A$, with the exception that
here we have to intersect with $\sp$.

We summarize the information we gathered about the Bruhat
decomposition for type $C$ in the following:

\begin{prop}
The group $\sp$ decomposes into double cosets of the form
$\mathbb{U}^C_{\pi} \pi \borelc$, where $\pi$ runs through $B_n$.
Every double coset decomposes into cosets of the form $A\borelc$
where $A$ is a general representative of the form
(\ref{representative}). The number of free parameters in $A$ is
equal to $n^2-\ell(\pi)$.
\end{prop}

\section{Sign Balance for Type $A$}\label{SignA}

Let $p$ be a prime number and let $q$ be a power of $p$. Denote by
$\mathbb{F}_q$ the field with $q$ elements. We prove the
following:

\bthm \label{mainA}
$$\sumlim_{K \in GL_n(\mathbb{F}_q)}{\omega_q^{s(K)}}=-(q-1)^{n-1}q^{n \choose 2}[n-1]_q!.$$
where $s(K)$ is the sum in $\mathbb{F}_q$ of the images of the
elements of the matrix $K$ under any bijection between
$\mathbb{F}_q$ and the set $\{0,..,q-1\}$ sending $0$ of
$\mathbb{F}_q$ to $0$ and $\omega_q$ is a primitive complex $q$-th
root of unity. \ethm

The following corollary is immediate:

\begin{cor}
The number of matrices in $GL_n(\mathbb{F}_q)$ whose sum of
entries is $0$ in $\mathbb{F}_q$ is exactly
$$[n-1]_q!(q-1)^{n-1}q^{{n \choose 2}}(q^{n-1}-1)$$ while for
every $1 \leq i \leq q-1$, the number of matrices in
$GL_n(\mathbb{F}_q)$ whose entries add up to $i$ in $\mathbb{F}_q$
is:
$$[n-1]_q!(q-1)^{n-1}q^{{n \choose 2}+n-1}.$$
\end{cor}\qed





Along this chapter we call a matrix or a column of a matrix
\emph{odd} if its sum of entries is not $0$ in $\mathbb{F}_q$ and
\emph{even} otherwise. A coset consisting entirely of odd matrices
will be called an \emph{odd coset}. In order to prove the theorem,
we take the following direction: Instead of summing over the whole
group of matrices, we sum over every coset separately. It turns
out that some of the cosets are sign-balanced in the sense that
they contain for each $i$ the same number of matrices with sum of
entries equals to $i$, while the others have only odd matrices.

 The following lemma identifies the sign-balanced cosets.

\blem \label{sign balanced cosets type A}

Let $A$ be a general representative of the double coset
$U_{\pi}\pi \mathbb{B}^+$ corresponding to $\pi \in S_n$.
 Make some substitution in the free parameters of $A$ to get a coset representative and call
it $\tilde{A}$. If $\tilde{A}$ has an odd column which is not the
last one then the coset $[\tilde{A}]=\{\tilde{A}B|B \in \borel\}$
is sign-balanced, i.e.,
$$\sumlim_{K \in \tilde{A}\mathbb{B}^+}{\omega_q^{s(K)}}=0.$$
\elem

\begin{proof}
Denote by $j$ the first odd column of $\tilde{A}$. For every
matrix $X \in \tilde{A}\borel$, the $k$-th column of $X$ is a
linear combination of the first $k$ columns of $\tilde{A}$. Now,
for every
 $B=B^{(0)} \in \borel$, construct the matrices $B^{(i)} \in \borel$ ($1 \leq i \leq q-1$),
 such that $B^{(i)}$ differs from $B$ only in the entry in location
 $(j,n)$ and $B^{(i)}_{j,n} \neq B^{(k)}_{j,n}$ if $i \neq k$.
 Note that the set $\{s(\tilde{A}B^{(i)})\mid 0 \leq i \leq q-1\}$ forms a full system of
representatives of $\mathbb{F}_q$. This gives us a partition of
the coset $\tilde{A} \borel$ into $q$ equal pieces, each of size
$\frac{q^{n \choose 2}(q-1)^n}{q}$. (Note that
$|\tilde{A}\borel|=|\borel|=(q-1)^n q^{{n \choose 2}}$).

\begin{exa}
The following example illustrates the bijection where
$G=GL_4(\mathbb{Z}_2)$.

 Let \begin{equation*} \pi=\begin{pmatrix}
 1 & 2 & 3 & 4 \\
 3 &1 & 2 &4
 \end{pmatrix} \in S_4.
 \end{equation*}

 Then a possible
representative of one of the cosets of $ \mathbb{U}_{\pi}\pi
\borel$ is:
  \begin{equation*}\tilde{A}=\begin{pmatrix}
  0 & 1 & 0 & 0 \\
  0 &  0 & 1 & 0 \\
  1 & 0 & 0 & 0\\
  1 & 0 & 1 & 1
  \end{pmatrix}.\end{equation*}
Column $j=2$ of $\tilde{A}$ is odd.

Now consider the Borel matrix:

  $$B=\begin{pmatrix}

               1 &  1 & 1 & 0 \\
               0 & 1 & 0 & 1\\
                0 & 0 & 1 & 0\\
                0 & 0 & 0 & 1

                    \end{pmatrix}$$

                    which produces the element
                    $$\tilde{A}B=\begin{pmatrix}
                    0 &  1 & 0 & 1 \\
                    0 & 0 & 1 & 0\\
                    1 & 1 & 1 & 0\\
                    1 & 1 & 0 & 1
                    \end{pmatrix}.$$

Toggle the element $B_{2,4}$ to get

                    $$B'=\begin{pmatrix}
                    1 &  1 & 1 & 0 \\
                    0 & 1 & 0 & 0\\
                    0 & 0 & 1 & 0\\
                    0 & 0 & 0 & 1
                    \end{pmatrix}$$

                    and the corresponding element of $\tilde{A}\mathbb{B}^+$ is

                     $$\tilde{A}B'=\begin{pmatrix}
                    0 &  1 & 0 & 0 \\
                    0 & 0 & 1 & 0\\
                    1 & 1 & 1 & 0\\
                    1 & 1 & 0 & 1
                    \end{pmatrix}.$$

Note that $\tilde{A}B$ and $\tilde{A}B'$ have opposite parity.

\end{exa}

Now Calculate:
\begin{align*}
\sumlim_{K \in \tilde{A}\mathbb{B}^+}{\omega_q^{s(K)}} & =  \sumlim_{i=0}^{q-1}{\sumlim_{\substack{K \in \tilde{A} \borel \\ s(K)=i}}{\omega_q^i}}\\
 & = (q-1)^n q^{{n \choose 2 }-1}\left(\sumlim_{i=0}^{q-1}{\omega_q^{i}}\right) = 0
\end{align*}
\end{proof}

\blem \label{contribution of odd cosets A}

Let $\pi \in S_n$. Let $A$ be a general representative of the
double coset $\mathbb{U}_{\pi}\pi \borel$ corresponding to $\pi
\in S_n$.
 Make some substitution in the free parameters of $A$ to get a coset representative, and call
it $\tilde{A}$. If the first $n-1$ column sums of $\tilde{A}$ are
$0 \pmod p$, then the imbalance calculated inside
$\tilde{A}\borel$ is:
$$\sumlim_{K \in \tilde{A}\borel}{\omega_q^{s(K)}}=-(q-1)^{n-1}q^{n \choose 2}.$$
\elem

\begin{proof}
The last column of $\tilde{A}$ contains only one nonzero element,
which must be $1$. (Indeed, by the construction of $A$, the last
column of $A$ contains a '$1$', coming from the permutation $\pi$,
say, in the row numbered $k$. The entries located in location
$(l,n)$ for $l<k$ are zeroes by the definition of $A$. On the
other hand, every row numbered $l>k$ must contain a '$1$' coming
from $\pi$ which is located in some column $j<n$. All of the
entries of the row $l$ to the right of $j$, including the entry
$(l,n)$  must be zero).\\
 Let $K=\tilde{A}B \in \tilde{A}\borel$.
 Since the first $n-1$ columns of $\tilde{A}$ add up to \mbox{$0 \pmod p$,}
$s(K)$ depends uniquely on the element of $B$ located in the entry
numbered $(n,n)$, which can be any element of $\mathbb{F}_q
\setminus \{0\}$. For every $0 \neq \alpha \in \mathbb{F}_q$ there
are $q^{n \choose 2}(q-1)^{n-1}$ Borel matrices $B=(b_{ij})$ with
$b_{nn}=\alpha$, and thus the imbalance of this coset is:
\begin{eqnarray}
\sumlim_{K \in \tilde{A}\borel}{\omega_q^{s(K)}} & = &
(q-1)^{n-1}q^{n \choose 2}(\omega_q+\omega_q^2+ \cdots
+\omega_q^{q-1}) \nonumber
\\ & = & -(q-1)^{n-1}q^{n \choose 2} \nonumber \qquad
\end{eqnarray}
\end{proof}

Until now we have seen that there are two types of cosets. A sign-
balanced coset and an odd coset. Only the last type contributes to
the imbalance $\sumlim_{K \in GL_n(\mathbb{F}_q)}{\omega_q^{s(K)}
}$. Since every coset contains $q^{{n \choose 2}}$ elements, it
remains to count the number of odd cosets. The following lemma
tells us which double cosets contain odd cosets.

\blem \label{identify odd cosets A} Let $\pi \in S_n$.  The double
coset $U_{\pi}\pi \mathbb{B}^+$ contains odd cosets if and only if
$\pi(n)=n$.  \elem

\begin{proof}

Let $A$ be a general representative of the double coset $U_{\pi}
\pi \borel$. If $\pi(n)=n$ then since the entries of the last row
of $A$ are not located above or to the right of any '$1$', the
last row of $A$ has no limitation on its parameters. Hence, every
column, except for the last one, has parameters which we can
choose such that it will be even. These choices yield odd cosets.

On the other hand, if $\pi(n) \neq n$ then there is some $i<n$
such that $\pi(i)=n$. This forces the $i$-th column of $A$ to be
odd and the corresponding coset to be sign-balanced.
\end{proof}

We turn now to the proof of Theorem \ref{mainA}: In order to
calculate the imbalance we have to count only non-balanced cosets.
By Lemma \ref{identify odd cosets A} we are interested only in the
double cosets corresponding to permutations $\pi \in S_{n-1}$.
Every such double coset has a total of ${n \choose 2}-\ell(\pi)$
parameters, which amounts to $q^{{n \choose 2}-\ell(\pi)}$
different cosets. In order to get an odd coset with sum of entries
not equal to $0$ modulo $q$ we have to choose the parameters in
such a way that every column except for the last one will sum up
to $0 \pmod p$. Exactly $\frac{1}{q}$ of the choices in each
column give an odd column and thus there are $q^{{n \choose
2}-\ell(\pi)-(n-1)} $ such cosets in each double coset. By Lemma
\ref{contribution of odd cosets A}, each coset contributes $-q^{n
\choose 2}(q-1)^{n-1}$ and we have in total:

\begin{align*}
\sumlim_{K \in GL_n(\mathbb{F}_q)}{\omega_q^{s(K)}} & =  \sumlim_{\substack{\pi \in S_n\\ \pi(n)=n}}-(q-1)^{n-1}q^{n \choose 2}q^{{n \choose 2}-(n-1)-\ell(\pi)}  \\
 & =  -(q-1)^{n-1}q^{n \choose 2}\sumlim_{\pi \in S_{n-1}}{q^{{n \choose 2} -\ell(\pi)-(n-1)}}  \\
 & =  -(q-1)^{n-1}q^{n \choose 2}\sumlim_{\pi \in S_{n-1}}{q^{{{n-1} \choose 2} - \ell(\pi)}}  \\
 & =  -(q-1)^{n-1}q^{n \choose 2}\sumlim_{\pi \in S_{n-1}}{q^{\ell(\pi)}}  \\
 & =  -(q-1)^{n-1}q^{n \choose 2}[n-1]_q! \quad\quad
\end{align*}

Note that the forth equality of the last calculation follows from
the obvious bijection inside $S_{n-1}$ given by multiplying by the
longest permutation. \qed

\section{Sign Balance for Type $C$}\label{SignC}

In this section we prove the following result: \bthm
\label{mainC2}
$$\sumlim_{K \in Sp_{2n}(\mathbb{Z}_2)}{(-1)^{o(K)}}=
-2^{n^2}\cdot[2]_2[4]_2 \cdots [2n-2]_2$$ where $o(K)$ is the
number of $1$'s in $K$.\ethm

The following corollary is immediate:
\begin{cor}
The number of even matrices in $Sp_{2n}(\mathbb{Z}_2)$ is exactly
$$2^{n^2-1}[2]_2 \cdots [2n-2]_2([2n]_2-1)$$ while the number of
odd matrices is $$2^{n^2-1}[2]_2 \cdots [2n-2]_2([2n]_2+1).
\quad\quad \qed$$
\end{cor}

In proving the theorem, we use the same strategy used for type A.
We sum over each coset separately and distinguish between odd and
sign-balanced cosets.  The following lemma identifies the
sign-balanced cosets.

\blem \label{sign balanced cosets for type C z2} Let $A$ be a
general representative of the double coset $U^C_{\pi}\pi \borelc$
corresponding to $\pi \in B_n$. Make some substitution in the free
parameters of $A$ to get a coset representative, and call it
$\tilde{A}$.  If $\tilde{A}$ has an odd column which is not the
last one, then the coset $[\tilde{A}]=\{\tilde{A}B\mid B \in
\borelc\}$ is sign-balanced, i.e.,
$$\sumlim_{K \in \tilde{A}\borelc}{(-1)^{o(K)}}=0.$$
\elem
\begin{proof}
An element of $\borelc$ is an invertible upper triangular matrix
which is also symplectic. If we take $b$ to be an upper triangular
matrix with a set of columns $\{v_1,...,v_{2n}\}$ then, as was
stated in Section \ref{Background type C}, forcing it to be
symplectic is equivalent to imposing the equations (note that we
are working over $\mathbb{Z}_2$):
$$B(v_i,v_j) =  \left\{ \begin{array}{cc}
{1}      &  {i+j=2n+1} \\
{0}          & {i+j\neq 2n+1}
\end{array}
\right.
$$
As is easy to check, the equations of the form $B(v_i,v_i)=0$ are
trivial over $\mathbb{Z}_2$. The equations of the form
$B(v_i,v_{2n+1-i})=1$ are also trivial. (Indeed,
$B(v_i,v_{2n+1-i})=\sumlim_{k=1}^{2n}{b_{k,i}\cdot
b_{2n+1-k,2n+1-i}}$ but since $b$ is upper triangular, over
$\mathbb{Z}_2$ we have $b_{ii} \cdot b_{2n+1-i,2n+1-i}=1$ and the
other summands vanish since for $k>i$ one has $b_{ki}=0$ and for
$k>2n+1-i$ one has $b_{2n+1-k,2n+1-i}=0$).

Now, the only nontrivial equations involving the parameters
appearing in the last column are the ones of the form:
$$B(v_i,v_{2n})=0 ,(2 \leq i \leq 2n-1)$$ and each such equation can be written in such
a way that the parameters of the last column are free while the
parameters of the first row depend on them. Explicitly, we write
the equation $B(v_i,v_{2n})=0$ as
$$b_{1i}=\sumlim_{k=2}^{2n}{b_{ki}\cdot b_{2n+1-k,2n}} .$$

Note that the elements of the last column of the matrix $b$ have
no appearance as a part of a linear combination in any place other
than the first row. This is justified by the fact that every
nontrivial equation, involving the first row, which we have not
treated yet must be of the form $B(v_i,v_j)=0$ for $1 \leq i <j
\leq 2n-1$. Thanks to the upper triangularity of $b$, the elements
laying in the first row vanish in these equations.

Let us look at the following example:
                     $$b=\begin{pmatrix}
                    1 &  b_{12} & b_{13} & b_{14} \\
                    0 & 1 & b_{23} & b_{24}\\
                    0 & 0 & 1 & b_{34}\\
                    0 & 0 & 0 & 1
                    \end{pmatrix}.$$
The only nontrivial equations involving the last column are:
$B(v_2,v_4)=0$ and $B(v_3,v_4)=0$.

These equations can be written as:
\begin{align*}
b_{12} & = b_{34} \\
b_{13} & = b_{24}+b_{23}\cdot b_{34}
\end{align*}
so after intersecting with $Sp_{2n}(\mathbb{Z}_2)$, the matrix $b$
looks like:
                     $$b=\begin{pmatrix}
                    1 &  b_{34} & b_{24}+b_{23}\cdot b_{34} & b_{14} \\
                    0 & 1 & b_{23} & b_{24}\\
                    0 & 0 & 1 & b_{34}\\
                    0 & 0 & 0 & 1
                    \end{pmatrix}.$$
The elements of the last column appear only in the first row and
in the equations of the form $B(v_i,v_j)=0$ the elements located
in the first row vanish.

Note that in this case all of the parameters outside the first row
are free. This doesn't hold in general. Nevertheless, as we have
proven, we can arrange the parameters such that the elements of
the last column reappear only in the first row.

Returning now to the proof, we have two cases:

\begin{itemize}
\item The first column of $\tilde{A}$ is odd.
In this case we can use the element located in place $(1,2n)$ to
construct a bijection between odd and even matrices inside the
coset $\tilde{A} \borelc$. This is done in the same way described
earlier for type A: Divide $\borelc$ into two disjoint subsets:
$$\borelc_0=\{T=(t_{i,j}) \in \borelc \mid t_{1,2n}=0\}$$
$$\borelc_1=\{T=(t_{i,j}) \in \borelc \mid t_{1,2n}=1\}.$$

For every matrix $X \in \tilde{A}\borelc$, the $k$-th column of
$X$ is a linear combination of the first $k$ columns of
$\tilde{A}$. Now, due to the fact that the parameter appearing in
the location $(1,2n)$ has no other appearance, for every $B \in
\borelc_0$ there is some $B' \in \borelc_1$ such that $B$ and $B'$
differ only in the entry numbered $(1,2n)$.

Note that $\tilde{A}B$ and $\tilde{A}B'$ are obtained from
$\tilde{A}$ by the same sequence of column operations except for
the first column which was used in producing $AB$ but was not used
in producing $\tilde{A}B'$. Hence $\tilde{A}B$ and $\tilde{A}B'$
have opposite parity. This gives us a bijection between the odd
and the even matrices of the coset $A\borelc$.

\item The first column of $\tilde{A}$ is even.
Denoting  by $j$ the number of the first odd column of
$\tilde{A}$, we use the element located in place $(j,2n)$ to
construct a bijection between the odd and even matrices inside
$\tilde{A} \borelc$ in the same way as in the previous case. Note
that since the element located in place $(j,2n)$ in the matrices
of the Borel subgroup can reappear only in the first row, it
affects only the first column of $\tilde{A}$, which is even.
\end{itemize}
\end{proof}

We turn now to treat the odd cosets. \blem \label{contribution of
odd cosets type C} Let $\pi \in B_n$. Let $A$ be a general
representative of the double coset $U_{\pi}\pi \borelc$
corresponding to $\pi \in B_n$.
 Make some substitution in the free parameters of $A$ to get a coset representative, and call
it $\tilde{A}$. If all of the first $2n-1$ columns of $\tilde{A}$
are even then all of the matrices belonging to the coset
$\tilde{A}\borelc$ are odd. The imbalance calculated inside this
coset is:
$$\sumlim_{K \in \tilde{A}\borelc}{(-1)^{o(K)}}=-|\borelc|=-2^{n^2}.$$
\elem
\begin{proof}
The last column of $\tilde{A}$ is always odd and thus since all
other columns of $\tilde{A}$ are even, $\tilde{A}$ itself is an
odd matrix and the same holds for $\tilde{A}B$ for every $B \in
\borelc$. The size of the coset $\tilde{A}\mathbb{B}^+$ is
$2^{n^2}$, and the result follows.
\end{proof}

\blem \label{identify odd cosets C}

Let $\pi \in B_n$. The double coset $U_{\pi}^C \pi \borelc$
contains odd cosets if and only if $\pi(2n)=2n$. \elem
\begin{proof}
Let $A$ be a general representative of the double coset $U_{\pi}^C
\pi \borelc$. Write $U=A \pi^{-1}$. Then $U \in
\mathbb{U}_{\pi}^C$ is a lower triangular matrix and since
$\pi(2n)=2n$ (which implies also $\pi(1)=1$), the first column as
well as the last row of $U$ contain $2n-1$ parameters. Note that
$U^T \in \borelc$ and thus by the considerations described in
Lemma \ref{sign balanced cosets for type C z2}, the parameters
appearing in the last column of $U^T$ can reappear only in the
first row of $U^T$. We conclude that the parameters of the last
row of $U$ can reappear only in the first column of $U$. Now, for
every column numbered $2 \leq k \leq 2n-1$ in $U$ and for every
choice of the first elements of the column numbered $k$, we are
free to choose the parameter located in the bottom of this column,
$(2n,k)$,  in such a way that the column will be even. The
parameter located in the place $(2n,1)$ has no other appearance
and thus we can choose all of the first $2n-1$ columns of $U$ to
be even. Getting back to the general representative $A$, since
$\pi(2n)=2n$, we have also $\pi(1)=1$ and thus $A$ and $U$ differ
only in the columns $1<k<2n$ so that the proof works also for $A$.

On the other hand, if $\pi(2n) \neq 2n$ then $\pi$ contains a
column numbered $k<2n$ which has only one nonzero element, located
in place $(2n,k)$. By the construction of the general
representative $A$, there are only zeros above the $1$ coming from
the permutation and thus this odd column appears also in $A$. By
the previous lemma, the coset $\{\tilde{A}B\mid B \in \borelc\}$
is sign-balanced.
\end{proof}

Now, we have to count the imbalance on the odd cosets. By Lemma
\ref{identify odd cosets C} we are interested only in the double
cosets corresponding to the permutations $\pi \in B_{n-1}$. The
following lemma shows how to count.

\blem \label{number of odd cosets type C}

Let $\pi \in B_n$ such that $\pi(n)=n$ . The double coset
$\mathbb{U}_{\pi}^C \pi \borelc$ contains exactly
$2^{(n-1)^2-\ell(\pi)}$ odd cosets. \elem

\begin{proof}
Let $A$ be representative of the double coset $\mathbb{U}_{\pi}^C
\pi \borelc$. As was shown in the previous lemma, the parity of a
each one of the first $2n-1$ columns of $A$ is determined by the
free parameter in its bottom. Since there are a total of
$n^2-\ell(\pi)$ free parameters and exactly $2n-1$ 'bottom
parameters', the number of substitutions of parameters giving all
of the $2n-1$ first columns even is $2^{n^2-\ell(\pi)-(2n-1)}$.
This is also the number of odd cosets in the double coset
$\mathbb{U}_{\pi}^C$.

\end{proof}

We turn now to the proof of Theorem \ref{mainC2}. In order to
calculate the imbalance we have to count only odd cosets. By Lemma
\ref{identify odd cosets C}, we are interested only in the double
cosets corresponding to permutations $\pi \in B_{n-1}$. By Lemma
\ref{number of odd cosets type C}, every such double coset
contains $2^{(n-1)^2-\ell(\pi)}$ odd cosets. By Lemma
\ref{contribution of odd cosets type C}, each odd coset
contributes $-2^{n^2}$ to the imbalance, and we have in total:

\begin{align*}
\sumlim_{K \in Sp_{2n}(\mathbb{Z}_2)}{(-1)^{o(K)}} & =
\sumlim_{\substack{\pi \in B_n \\ \pi(n)=n}}{-2^{n^2} \cdot 2^{(n-1)^2-\ell(\pi)}}  \\
& =  -2^{n^2} \sumlim_{\pi \in B_{n-1}}{2^{(n-1)^2-\ell(\pi)}} \\
& =  -2^{n^2} \sumlim_{\pi \in B_{n-1}}{2^{\ell(\pi)}}  \\
& =  -2^{n^2}[n-1]_2!  \\
& =  -2^{n^2}\cdot[2]_2[4]_2 \cdots [2n-2]_2
\end{align*}

\section{Appendix}
Our results can be seen as an example of the {\bf Cyclic sieving
phenomenon} introduced in the paper of Reiner, Stanton and White
\cite{RSW}. In this section we present this point of view. We
start with the definition of the cyclic sieving phenomenon,
following \cite{RSW}.

Let $X$ be a finite set and let $C_n$ be the cyclic group acting
on $X$. Let $X(t)$ be a polynomial in $t$ having nonnegative
integer coefficients, with the property that $X(1)=|X|$. One can
think of $X(t)$ as a generating function for $X$. Fix an
isomorphism $\omega$ of $C_n$ with the complex $n$- th roots of
unity. Then it is easy to see that the following are equivalent:

\begin{enumerate}
\item For every $c \in C_n$:
$$[X(t)]_{t=\omega(c)}=|\{x \in X \mid c(x)=x\}|.$$

\item The coefficient $a_l$ defined uniquely by the expansion
$$X(t)= \sumlim_{l=0}^{n-1}{a_l t^l} \quad mod \quad t^{n} -1$$ has the
following interpretation: $a_l$ counts the number of $C_n$-orbits
on $X$ for which the stabilizer-order divides $l$. In particular,
$a_0$ counts the total number of $C_n$-orbits on $X$, and $a_1$
counts the number of free $C_n$ orbits on $X$.
\end{enumerate}

When either of these two conditions holds, one says that
$(X,X(t),C_n)$ has the {\bf Cyclic sieving phenomenon}.\\

In our case, $X$ is either the set $GL_n(\mathbb{F}_q)$ or
$Sp_{2n}(\mathbb{Z}_2)$ and $$X(t)=\sumlim_{K \in G}{t^{o(K)}}$$
where $G$ is one of the above groups while the group acting on $X$
is $C_q$ in the first case and $C_2$ in the second. We describe
now the action of $C_q$ in the type A case. The other case is very
similar.

Note first that if $q$ is not a prime number then we have to
define a linear order on $\mathbb{F}_q$. For example, consider
$\mathbb{F}_q$ as a vector space over $\mathbb{Z}_p$, $p$ prime
and take the lexicographic order of the coordinates).\\

Now, for every matrix $K \in GL_n(\mathbb{F}_q)$, we have two
cases:

\begin{itemize}

\item
$K=\tilde{A}B$ where $\tilde{A}$ is a representative of an odd
coset. In this case the action of $C_q$ is trivial.
\item $K=\tilde{A}B$ where $\tilde{A}$ is a representative of a sign-balanced
coset. In this case we denote by $j$ the first odd column of
$\tilde{A}$ which is not the first one and define for the
generator $c$ of $C_q$:
$$ c \cdot \tilde{A}B=AB'$$
where $B'$ is obtained from $B$ by replacing $B_{j,n}$ by the
successor of this element with respect to the prescribed order.

\end{itemize}

{\bf Acknowledgements:} We would like to thank Prof. Yuval Roichman and Prof. Ron Adin for
the time they spent in discussions on this work. We also thank Mishael Sklarz for
helping us in implementing the computerized side of this work.

 \end{document}